\documentclass[12pt]{amsart}
\usepackage[utf8]{inputenc}
\usepackage{setspace}
\usepackage[english]{babel}
\usepackage{mathabx}
\usepackage{framed}
\usepackage{soul}
\usepackage[all]{xy}
\usepackage[margin=1in]{geometry}
\usepackage{arydshln}
\usepackage[colorinlistoftodos]{todonotes}
\usepackage{pb-diagram}
\usepackage[mathscr]{euscript}
\usepackage{multicol}
\usepackage{stmaryrd} 
\usepackage{empheq} 
\usepackage{tikz-cd}
\tikzset{
  symbol/.style={
    draw=none,
    every to/.append style={
      edge node={node [sloped, allow upside down, auto=false]{$#1$}}}
  }
}
\usepackage{booktabs}

\usepackage{amsmath,amsfonts,amssymb,amsthm,mathrsfs,latexsym,mathtools}
\usepackage{commath}
\usepackage[T1]{fontenc}
\usepackage{hyperref}
\usepackage{caption} 
\usepackage{subcaption} 
\usepackage{booktabs,longtable,array}

\DeclareMathAlphabet{\mathbbmsl}{U}{bbm}{m}{sl}

\usetikzlibrary{matrix}

\newcommand{\C}{\mathbb{C}} 
\newcommand{\Z}{\mathbb{Z}}

\newcommand{\PP} {\mathbb{P}}

\newcommand{\A}{\mathcal A}
\newcommand{\M}{\mathcal M}
\newcommand{\J}{\mathcal J}
\newcommand{\Prym}{\operatorname{Prym}}

\newcommand{\Bl}{\operatorname{Bl}}

\newcommand{\ct}{\mathrm{ct}}

\usepackage{amscd} 

\usepackage{enumitem} 
\usepackage{comment}

\newtheorem{theorem}{Theorem}[section]
\newtheorem{definition}[theorem]{Definition}
\newtheorem{proposition}[theorem]{Proposition}

\newtheorem{question}[theorem]{Question}

\newtheorem{lemma}[theorem]{Lemma}

\newtheorem{remark}[theorem]{Remark}

\newtheorem{thm}{Theorem}

\tikzset{commutative diagrams/.cd,
mysymbol/.style={start anchor=center,end anchor=center,draw=none}
}

\title{Rationality of Fano threefolds and intermediate Jacobians}
\author{Yilong Zhang}
\address{Department of Mathematics,
Boyd Research and Education Center,
University of Georgia,
Athens, GA 30602, U.S.A.}
\address{Institute for Computational and Experimental Research
in Mathematics (ICERM), Brown University,
Providence, RI 02903, U.S.A.}
\email{yilong.zhang@uga.edu}

\date{Nov 17, 2026}

\begin{document}
\maketitle
\begin{abstract}
    In this paper, we survey the Clemens--Griffiths criterion for rationality of smooth Fano threefolds. We show this criterion is also sufficient for a general Fano 3-fold. In other words, a general Fano threefold is rational if and only if its intermediate Jacobian is isomorphic, as a principally polarized abelian variety, to a product of Jacobians of smooth curves. This theorem is based on results of Clemens--Griffiths, Beauville, Tyurin and Shokurov.
\end{abstract}

\section{Introduction}

Let $X$ be a smooth projective
threefold over $\C$ with $H^{1,0}(X)=H^{3,0}(X)=0$.  Its intermediate Jacobian
\[
 J(X)=H^3(X,\C)\big/(F^2H^3(X,\C)+H^3(X,\Z))
\]
is an abelian variety with a 
principal polarization induced by the cup-product pairing.

Throughout the paper, a \textit{product of Jacobians} means a product of Jacobians of smooth connected projective curves, with the product principal polarization. All isomorphisms of intermediate Jacobians are understood to preserve their principal polarizations. 

Clemens and Griffiths introduced a necessary condition for rationality.

\begin{theorem}[Clemens--Griffiths, \cite{CG}]
\label{Thm_CG_Criterion}
Let $X$ be a smooth projective rational threefold.  Then there are smooth
projective curves $C_1,\ldots,C_r$ such that
\[
 (J(X),\Theta_X)\simeq
 \prod_{i=1}^r (J(C_i),\Theta_{C_i})
\]
as principally polarized abelian varieties.
\end{theorem}

They used this criterion to prove that every
smooth cubic threefold is irrational \cite{CG}. We call the condition that $J(X)$ is not a product of Jacobians the \textit{intermediate Jacobian obstruction}. This suggests the following
question.

\begin{question}  \label{Q_CGFano}
Is the Clemens--Griffiths criterion sufficient for smooth Fano threefolds? In other words, if $X$ is a Fano 3-fold and $J(X)$ is a product of Jacobians, is $X$ rational?
\end{question}




There are 105 deformation families of smooth Fano 3-folds, and rationality is known (at least for a general member in a deformation class). Historically, several families of Fano 3-folds were shown to be irrational using other methods.  For example, Iskovskikh-Manin \cite{IM} showed the irrationality of quartic 3-folds using birational rigidity. Recently, Farb \cite{Far25} used the automorphism group of Klein quartic threefold to show its intermediate Jacobian is not a product of Jacobians. This shows Question \ref{Q_CGFano} holds at least for a general quartic threefold. 
This motivates the author to explore other families of Fano 3-folds.



We show that Question \ref{Q_CGFano} holds true for a general Fano 3-fold. Here \textit{general} means belonging to a nonempty Zariski-open subset of the deformation family.

\begin{theorem}
\label{thm:main}
Let $\mathfrak F$ be a deformation family of smooth complex Fano
threefolds, and let $X$ be a general member of $\mathfrak F$.  Then
\[
 X\text{ is rational}
 \quad\Longleftrightarrow\quad
 J(X)\text{ is a product of Jacobians of smooth curves}.
\]
\end{theorem}

This proof is a case-by-case argument based on \cite{CG}, \cite{Bea77}, \cite{Tyurin}, \cite{Sho83}. The main tools are the rationality criterion for conic bundles developed by Beauville and Shokurov \cite{Bea77,Sho83}, and degeneration of intermediate Jacobians, developed by Clemens and Tyurin \cite{Cle77,Tyurin}.


\subsection{Additional discussions}
\begin{remark}\normalfont
 One should note that the answer for Question \ref{Q_CGFano} can be no for a singular $X$, or weak Fano 3-folds. Voisin \cite{Voisin_Quartic7Nodes} showed a very general quartic double solid with 7 nodes is irrational (actually stably irrational), but has intermediate Jacobian in $\mathcal A_3$, which consists of product of Jacobians.
\end{remark}

\begin{remark}\normalfont

Theorem~\ref{thm:main} concerns the general member of each deformation
family.  In several cases, Question \ref{Q_CGFano} is known for all smooth members, like cubic threefolds \cite{CG}, quartic double
solids \cite{Voi88}. Both require a deeper understanding of theta divisors for every smooth member. 
Question \ref{Q_CGFano} still remains open for all smooth members for the following cases: sextic double solid, quartic 3-fold, (2,3) complete intersection in $\mathbb P^5$; Gushel-Mukai 3-fold (note that the rationality question is not completely settled for special smooth members \cite{Debarre26}), double Veronese cone. See Table \ref{tab:rho1}.
\end{remark}

\begin{remark}\normalfont
By \cite{HT19} and \cite{Engel2025Matroids}, a very general irrational Fano 3-fold is also stably irrational, so Theorem \ref{thm:main} also shows that for a very general Fano 3-fold $X$,
$$X\ \textup{is stably irrational}\Leftrightarrow J(X)\ \textup{is not a product of Jacobians of curves}.$$
\end{remark}

\noindent\textbf{Structure of the paper}: In Section \ref{sec:setup}, we review basic facts about period images for curves of compact type and the closure in Satake--Baily--Borel compactification; We review degeneration techniques, and Beauville's criterion for conic bundles. In Section \ref{sec:cases}, we summarize the known results to Question \ref{Q_CGFano} and prove Theorem \ref{thm:main}. In Section \ref{sec:sextic}, we demonstrate the proof in the sextic double solid case. \\

\noindent \textbf{Acknowledgements}. I would like to thank Nick Addington, Valery Alexeev, Benson Farb, and Lena Ji for useful conversations. I would like to thank ICERM for its hospitality during the semester program \textit{Computations on K3 Surfaces and Related Varieties}. \\

\noindent \textbf{AI Declaration}. The author used ChatGPT 5.6 Sol to assist in working out the proof of Lemma \ref{lemma_sextic_double_solid_degeneration}, checking references, and improving the English. The manuscript was drafted by the author. The author independently verified all mathematical statements and references and takes full responsibility for the content of the paper.

\section{Preliminaries}
\label{sec:setup}

\subsection{The (compact-type) Torelli locus}

Let $\A_g$ be the moduli space of principally polarized abelian varieties of
dimension $g$.  For $g\ge 2$, let $\mathcal M_{g}^{ct}\subset \overline{\mathcal M}_g$ be the locus of stable curves of compact type, that is, stable curves whose dual graph is a tree. For such a curve $C$, its generalized
Jacobian is an abelian variety and
\[
 J(C)=\prod_{v}J(C_v),
\]
where $C_v$ are the irreducible components of $C$.  The extended Torelli map
therefore defines
\[
 \tau_g^{\ct}:\overline{\M}_g^{\ct}\longrightarrow\A_g.
\]

\begin{definition}\normalfont
The \emph{compact-type Torelli locus} is
\[
 \J_g^{\ct}:=\tau_g^{\ct}(\overline{\M}_g^{\ct})\subset\A_g.
\]
Equivalently, it is the locus of product of Jacobians of total dimension $g$.

\end{definition}

Recall that the Satake--Baily--Borel (SBB) compactification $\overline{\A}_g^{SBB}$ of $\mathcal A_g$ is  
\[
 \overline{\A}_g^{SBB}=\A_g\sqcup\A_{g-1}\sqcup\cdots\sqcup\A_0
\]
 and the boundary $\partial \overline{\A}_g^{SBB}=\A_{g-1}\sqcup\cdots\sqcup\A_0$ consists of lower dimensional $\mathcal A_{g'}$. The boundary contains the pure weight-one part of a limiting mixed Hodge structure. 

\begin{proposition}[{\cite{Hoyt}, \cite[Proposition 3.2.2]{Tyurin}}]
\label{prop:satake-closure}
The locus $\J_g^{\ct}$ is closed in $\A_g$. Its closure in the Satake--Baily--Borel compactification $\overline{\A}_g^{SBB}$ is
\[
 \overline{\J_g^{\ct}}
 =\bigsqcup_{h=0}^{g}\J_h^{\ct}.
\]
In particular, a SBB limit of products of Jacobians is again a product of Jacobians, possibly of smaller dimension. 
\end{proposition}
\subsection{Griffiths component} In this subsection, we give a modern overview of Clemens--Griffiths criterion (cf. Theorem \ref{Thm_CG_Criterion}). $X$ is a smooth projective threefold with $H^{1,0}(X)=H^{3,0}(X)=0$. 
Every principally polarized abelian variety decomposes uniquely, up to the
order of the factors, into indecomposable principally polarized abelian
varieties.  The
\emph{Griffiths component} $J_G(X)$ is the product of those indecomposable
factors of $J(X)$ which are not Jacobians of smooth curves \cite{Tyurin}.  

\begin{proposition}[{\cite[Theorem 3.25]{CG}}]
\label{prop:blowup}

$J_G(X)$ is a
birational invariant of smooth projective threefolds.

\end{proposition}

\begin{proof}
Let $Y$ be a smooth projective threefold and let $\widetilde Y=\Bl_ZY$. The blow-up formula gives
\[
 H^3(\widetilde Y,\Z)\simeq
 H^3(Y,\Z)\oplus H^1(Z,\Z)(-1)
\]
when $Z$ is a curve, and no new degree-three summand when $Z$ is a point.
The decomposition respects the Hodge structures and the intersection
pairings, giving the asserted product decomposition with the product
principal polarization.  A birational map between smooth projective
threefolds can be factored into blow-ups and blow-downs with smooth centers
\cite{AKMW02}\footnote{Of course, the original proof in \cite{CG} only uses Hironaka's resolution of singularities.}.
Since a blow-up only adds a curve Jacobian and the indecomposable
decomposition of a principally polarized abelian variety is unique, the
non-Jacobian factors are unchanged.  
\end{proof}
Let $g=\dim J(X)$, we obtain
\[
 J_G(X)=0\quad\Longleftrightarrow\quad [J(X)]\in\J_g^{\ct}.
\]
Since $J_G(\mathbb P^3)=0$, this implies Theorem \ref{Thm_CG_Criterion}.

\subsection{Degeneration criterion} 

Suppose now that $\mathcal X^*\to\Delta^*$ is a smooth projective family of
threefolds with $H^{1,0}(X_t)=H^{3,0}(X_t)=0$, and let $g=h^{2,1}(X_t)$.  After a finite base change its monodromy is
unipotent, and the family of intermediate Jacobians has semiabelian
reduction
\[
 1\longrightarrow(\C^*)^r\longrightarrow J_{\lim}
 \longrightarrow A_{\lim}\longrightarrow 0.
\]
The point $[A_{\lim}]\in\A_{g-r}\subseteq \overline{\A}_g^{SBB}$ is the SBB limit.  The abelian variety
$A_{\lim}$ is called the \emph{compact part} of the limiting intermediate Jacobian.

\begin{proposition}[Degeneration criterion]
\label{prop:degeneration}
In the preceding situation, assume that the compact part
\begin{equation}\label{eqn_compact_lim}
     [A_{\lim}]\notin\J_{g-r}^{\ct}.
\end{equation}

Then $[J(X_t)]\notin\J_g^{\ct}$ for general $t\in\Delta^*$.  Consequently
$X_t$ is irrational for general $t$.
\end{proposition}

\begin{proof}
If $J(X_t)$ belonged to $\J_g^{\ct}$ for general $t$, its SBB limit would
belong to $\J_{g-r}^{\ct}$ by Proposition~\ref{prop:satake-closure}, contrary
to the hypothesis.  The final assertion follows from
Theorem~\ref{Thm_CG_Criterion}.
\end{proof}

\begin{remark}\normalfont
  Note that the degeneration criterion Proposition \ref{prop:degeneration} only uses the data in compact type, equivalently, the graded weight-one piece $gr_1^WH^3(X_t)_{lim}(1)$ of the limiting mixed Hodge structure, and does not use the extension data. This is enough for remaining cases to prove Theorem \ref{thm:main}. However, we want to point it out that Proposition \ref{prop:degeneration} is not enough to show a general cubic threefold is irrational. This is because a nodal cubic threefold is a Jacobian of a genus four curve, so \eqref{eqn_compact_lim} is not satisfied.
  
   There are other degeneration techniques devoted to giving a new proof of the irrationality of a general cubic threefold. Collino \cite{Collino79_irrationality}  and Bardelli \cite{Bar84} used extension data for the limiting intermediate Jacobian. Gwena \cite{Gwena} studied the combinatorial data of monodromy associated to a degeneration to the Segre cubic threefold with ten nodes. The technique is further developed in \cite{Engel2025Matroids} to show stable irrationality of a very general cubic threefold. 
\end{remark}

\subsection{Conic bundles and Prym varieties}

Let $q:Y\to S$ be a standard conic bundle over a smooth rational surface.
Let $D\subset S$ be its discriminant curve, and let
$\widetilde D\to D$ be the double cover which records a component of a
reducible conic.  Under the standard smoothness assumptions this is an
admissible double cover, and
\[
 J(Y)\simeq\Prym(\widetilde D/D)
\]
as principally polarized abelian varieties \cite{Bea77}.  We shall use the
following theorem.

\begin{theorem}[\cite{Bea77}]
\label{thm:beauville-prym}Let $q:Y\to \mathbb P^2$ be a standard conic bundle with discriminant $D\subset \mathbb P^2$ of degree at least six, then
$$ J(Y)\simeq\Prym(\widetilde D/D)$$
is not a product of Jacobians.
In particular,  $Y$ is irrational.
\end{theorem}

 Shokurov showed the following converse to the Clemens–Griffiths criterion in this setting.
\begin{theorem}[{\cite[10.1]{Sho83}}, {\cite[Theorem 9.1, 9.2]{Pro18}}]\label{thm_Sho}
    Let $\pi:X\to S$ be a standard conic bundle over a minimal rational surface. Then $X$ is rational if and only if $J(X)$ is a product of Jacobians. 
\end{theorem}

\section{The Fano classification and the known cases}
\label{sec:cases}

The classification of smooth Fano threefolds is due to
Iskovskikh in Picard rank one and to Mori--Mukai in higher Picard rank; see
\cite{IP99,MM81} and the online tables \cite{Fanography}.  The rationality
organization used below follows \cite[Section~2]{HT19}.

Note that, Question \ref{Q_CGFano} holds true simultaneously in a given birational class, according to Proposition \ref{prop:blowup}. Therefore we only need to consider minimal Fano 3-folds. 

\subsection{Picard rank one}

Table~\ref{tab:rho1} lists the minimal irrational families of Picard rank
one.  Here $h^{2,1}(X)$ is the dimension of the intermediate Jacobian $J(X)$. The labels come from \cite{Fanography}. The geometric models in the table describe general members where indicated.  In the last column, we summarize current status of Question \ref{Q_CGFano} for general member or all smooth members, give reference, and the techniques. 

Note that degeneration and automorphism techniques only solve Question \ref{Q_CGFano} for a general member. Settling the all smooth member relies on understanding of singularities of theta divisors.

\begin{longtable}{@{}>{\raggedright\arraybackslash}p{0.07\textwidth}
                        >{\raggedright\arraybackslash}p{0.30\textwidth}
                        >{\centering\arraybackslash}p{0.08\textwidth}
                        >{\centering\arraybackslash}p{0.17\textwidth}
                        >{\raggedright\arraybackslash}p{0.25\textwidth}@{}}
\caption{Minimal irrational Fano threefolds of Picard rank one.}
\label{tab:rho1}\\
\toprule
Label & Family & $h^{2,1}$ & method & IJ obstruction\\
\midrule
\endfirsthead
\toprule
Label & Family & $h^{2,1}$ & method & IJ obstruction\\
\midrule
\endhead
\bottomrule
\endfoot
$1$--$1$ & Sextic double solid
& $52$ & degeneration &general \cite{Tyurin}\\ 
$1$--$2$ & Quartic hypersurface in $\PP^4$
& $30$ & degeneration; automorphism &general \cite{Tyurin,Bar84,Far25}\\
$1$--$3$ & Complete intersection $X_{2,3}\subset\PP^5$
& $20$ & degeneration &general \cite{Bea77}\\
$1$--$4$ & Complete intersection of three quadrics in $\PP^6$
& $14$ & conic bundle&smooth \cite{Bea77}\\
$1$--$5$ & Prime Fano threefold of genus $6$ (ordinary Gushel--Mukai
threefold)
& $10$ & degeneration&general \cite{Bea77}\\
$1$--$7$ & Prime Fano threefold $V_{14}$ of genus $8$
& $5$ & birational to a cubic&smooth \cite{Tak89,IM00}\\
$1$--$11$ & Double Veronese cone, equivalently a smooth sextic in
$\PP(1,1,1,2,3)$
& $21$ & degeneration &general \cite{Tyurin} 
\\
$1$--$12$ & Quartic double solid
& $10$ & theta divisor &smooth \cite{Bea77,Voi88}\\
$1$--$13$ & Cubic threefold
& $5$ & theta divisor &smooth \cite{CG,Bea81}\\
\end{longtable}

For a smooth cubic threefold $X$ (family 1-13), Clemens and Griffiths \cite{CG} studied the Gauss map of the theta divisor $\Theta_X$ using Abel-Jacobi maps and computed its degree to be 72. This differs from theta divisors of a genus 5 Jacobian. Later, Beauville \cite{Bea81} showed that the singular loci $\Theta_X^{sing}$ is a unique point $p$, of multiplicity three, and that the projective tangent cone $\mathbb P T_p\Theta_X$ is isomorphic to the cubic threefold $X$ itself. However, singular locus of theta divisor of a genus 5 Jacobian has dimension at least one. This gives a second proof of sufficiency of Clemens-Griffiths criterion, and therefore the irrationality of cubic threefolds.

For quartic double solid case (family 1-12), Beauville \cite{Bea77} first proved the intermediate Jacobian obstruction for a general member using degeneration. Then later, Voisin \cite{Voi88} studied the singularities of the theta divisor for every smooth member. She showed that the singular locus has an irreducible component of dimension $5=h^{2,1}-5$. This excludes product of Jacobians and settles Question \ref{Q_CGFano} for all smooth quartic double solids.

\subsection{Higher Picard rank}

The minimal irrational families of Picard rank at least two are listed in
Table~\ref{tab:rhohigher}.  Question \ref{Q_CGFano} is known to be true for all smooth members. This is because all of these Fano 3-folds admit standard conic-bundle models over
rational surfaces \cite{IP99,HT19}. Then the theorem for these families of Fano threefolds follows from Theorem \ref{thm_Sho}. The first three cases also follow from Theorem \ref{thm:beauville-prym}.

\begin{longtable}{@{}>{\raggedright\arraybackslash}p{0.09\textwidth}
                        >{\raggedright\arraybackslash}p{0.44\textwidth}
                        >{\centering\arraybackslash}p{0.05\textwidth}
                        >{\raggedright\arraybackslash}p{0.25\textwidth}@{}}
\caption{Minimal irrational Fano threefolds of Picard rank at least two.}
\label{tab:rhohigher}\\
\toprule
Label & Family & $h^{2,1}$ & Discriminant\\
\midrule
\endfirsthead
\toprule
Label & Family & $h^{2,1}$ & Conic bundle discriminant\\
\midrule
\endhead
\bottomrule
\endfoot
$2$--$2$ & Double cover of $\PP^1\times\PP^2$ branched along a divisor of
bidegree $(2,4)$ & $20$ & Degree 8 in $\mathbb P^2$\\
$2$--$6$ & Verra threefold: a divisor of bidegree $(2,2)$ in
$\PP^2\times\PP^2$ & $9$ & Degree 6 in $\mathbb P^2$\\
$2$--$8$ & Double cover of $\Bl_p\PP^3$ branched along an anticanonical
divisor & $9$ & Degree 6 in $\mathbb P^2$\\
$3$--$1$ & Double cover of $(\PP^1)^3$ branched along a divisor of
multidegree $(2,2,2)$ & $8$ & Bidegree $(4,4)$ in $\mathbb P^1\times \mathbb P^1$\\
\end{longtable}

For example, projection of a Verra threefold to one factor of $\PP^2$ is a conic bundle with sextic discriminant. Family $2$--$6$ also contains double covers of the flag threefold $F(1,2)\subset\PP^2\times\PP^2$, branched along an anticanonical divisor. These also have conic-bundle structures over $\PP^2$.

\subsection{Proof of the main theorem}
Now we prove the main theorem (Theorem \ref{thm:main}).
\begin{theorem}
    A general Fano threefold $X$ is irrational if and only if $J(X)\notin \mathcal J_g^{ct}$.
\end{theorem}
\begin{proof}
    The "if" direction follows from Theorem \ref{Thm_CG_Criterion}. For "only if" direction, by Proposition \ref{prop:blowup}, it suffices to consider minimal irrational Fano 3-fold $X$. By the classification, it reduces to the cases in Table \ref{tab:rho1} and \ref{tab:rhohigher}.

The four families in Table \ref{tab:rhohigher} admit standard conic-bundle models over minimal rational surfaces, so Theorem \ref{thm_Sho} applies. This also includes family 1--4 in Table \ref{tab:rho1}. 

Clemens--Griffiths treat all smooth cubic threefolds (family 1-13). The results for $V_{14}$ (family 1–7) follows from its birational correspondence with a smooth cubic threefold \cite{Tak89,IM00}. Beauville’s nodal degenerations treat general (2,3) complete intersections (family 1–3), general Gushel–Mukai threefolds (family 1–5), and general quartic double solids (family 1–12) \cite[Theorem 5.6]{Bea77}. Voisin’s theorem gives the stronger assertion for every smooth
quartic double solid \cite{Voi88}.
    
The degeneration results for general sextic double solids (family 1–1), quartic threefolds
(family 1–2), and double Veronese cones (family 1–11) are discussed by Tyurin \cite[p.738-740]{Tyurin}. See also \cite{Bar84}, \cite[p.78]{ChePar}, and \cite[Section 7]{Pro26}. Farb \cite{Far25} gives a different proof for general quartic threefolds. This exhausts the families in the two tables and proves the theorem.
\end{proof}
\section{Sextic double solid}
\label{sec:sextic}

Tyurin sketched a degeneration argument for the intermediate Jacobian obstruction for general sextic double solid (and also for quartic threefolds, double Veronese cones, etc) \cite[p.738]{Tyurin}, but does not contain many details. In this section, we demonstrate the proof in the sextic double solid case, and work out the degeneration explicitly.



\begin{lemma}[{\cite{Tyurin}}] \label{lemma_sextic_double_solid_degeneration} There exists a sextic double solid $X_0$ whose branch surface has an ordinary quadruple point and whose resolution $\tilde{X}_0$ is a conic bundle over $\mathbb P^2$ with smooth discriminant of degree ten. In particular, $J(\tilde{X}_0)$ is not a  product of Jacobians.
\end{lemma}
\begin{proof}
     We consider degeneration of the sextic surface to a singular sextic $B_0$
$$x_3^2f_4(x_0,x_1,x_2)+x_3f_5(x_0,x_1,x_2)+f_6(x_0,x_1,x_2)=0.$$
The singularity is at $p=[0:0:0:1]$ which is an ordinary quadruple point. We can blow up the singularity and resolve the singularity on the double cover $X_0$ and get $\tilde{X}_0$. 

Explicitly, let
\[
 X_0=(w^2=x_3^2f_4+x_3f_5+f_6)
 \subset\PP(1,1,1,1,3)
\]
be the corresponding double solid, and let $q\in X_0$ be the point over
$p$.  In local coordinates at $q$ its equation is
\begin{equation}
\label{eq:local-sing}
 w^2=f_4(x_0,x_1,x_2)+f_5(x_0,x_1,x_2)+f_6(x_0,x_1,x_2).
\end{equation}

Give $(x_0,x_1,x_2,w)$ weights $(1,1,1,2)$ and perform the corresponding
weighted blow-up.  Write locally $x_i=s y_i$ and $w=s^2u$.  The strict
transform has equation
\begin{equation}
\label{eq:weighted-blowup}
 u^2=f_4(y_0,y_1,y_2)+s f_5(y_0,y_1,y_2)
                     +s^2f_6(y_0,y_1,y_2).
\end{equation}
Its exceptional divisor is
\[
 E=(u^2=f_4(y_0,y_1,y_2))\subset\PP(1,1,1,2),
\]
a smooth del Pezzo surface of degree two.  Since $E$ is smooth, equation
\eqref{eq:weighted-blowup} shows that the strict transform is smooth along
$E$; away from $E$ it is isomorphic to $X_0\setminus\{q\}$.  We have
therefore obtained a resolution
\[
 \rho:\widetilde X_0\longrightarrow X_0
\]
with rational exceptional divisor $E$.

Note that we can project $B_0$ from $p$, this induces a conic bundle structure on $\pi:\tilde{X}_0\to \mathbb P^2$, which has discriminant locus being a planar curve $D$ of degree ten:
$$f_5^2-4f_4f_6=0.$$

For general $f_4,f_5,f_6$, $D$ is smooth (Bertini theorem), and the fibers over a point in $D$ is union of two distinct lines. Marking the irreducible components in the singular conic defines an \'etale double cover 
$$\tilde{D}\to D.$$
Beauville's conic-bundle calculation now gives an isomorphism of principally
polarized abelian varieties
 $$J(\widetilde X_0)\simeq\Prym(\widetilde D/D).$$

Direct computation shows that $g(D)=(10-1)(10-2)/2=36$ and $g(\tilde{D})=71$. Therefore, the Prym variety $Prym(\tilde{D}/D)$ has dimension $35$. According to Beauville (Theorem \ref{thm:beauville-prym}), this Prym variety is not a product of Jacobians.
\end{proof}

\begin{proposition}[\cite{Tyurin}]\label{prop_sextic_double_solid}
    The intermediate Jacobian of a general smooth sextic double solid is not a product of Jacobians. 
\end{proposition}
\begin{proof}

 After performing a semistable reduction of the family in Lemma \ref{lemma_sextic_double_solid_degeneration}, the compact part of the limiting intermediate Jacobian has a direct summand $P=Prym(\tilde{D}/D)$. The above shows that $P\notin \mathcal{J}_{35}^{ct}$ is not in the compact-type Torelli locus. Hence the compact part of the limiting intermediate Jacobian, i.e., the limit point in the SBB boundary, is not in the $\mathcal J^{ct}$. This shows that the intermediate Jacobian $J(X)$ of a general sextic double solid $X$ is not a product of Jacobians.
\end{proof}

\begin{remark}\normalfont
    In fact, a semistable model for this degeneration has central fiber
    $$\tilde{X}_0\cup_EV,$$
    where $\tilde{X}_0$ is the desingularized sextic double solid, $V$ is a quartic double solid, and they intersect along $E$, which is a del Pezzo surface of degree two.

So the compact part of the limiting intermediate Jacobian is the product ppav
$$A_{\lim}\cong Prym(\tilde{D}/D)\times J(V),$$
of dimension $35+10=45$, and the semiabelian limit has torus rank 7. Neither factor is a product of Jacobians of curves, so $A_{\lim}\notin\mathcal J^{ct}$.
\end{remark}

\bibliographystyle{plain}
\bibliography{bibfile}

@incollection {Bea16,
    AUTHOR = {Beauville, Arnaud},
     TITLE = {The {L}\"{u}roth problem},
 BOOKTITLE = {Rationality problems in algebraic geometry},
    SERIES = {Lecture Notes in Math.},
    VOLUME = {2172},
     PAGES = {1--27},
 PUBLISHER = {Springer, Cham},
      YEAR = {2016},
   MRCLASS = {14M20 (14E08)},
  MRNUMBER = {3618664},
MRREVIEWER = {Alexander Duncan},
}

@article{Debarre26,
  author  = {Debarre, Olivier},
  title   = {{Gushel--Mukai} varieties},
  journal = {EMS Surv. Math. Sci.},
  year    = {2026},
  note    = {Published online first},
  doi     = {10.4171/EMSS/125}
}

@article {IM,
    AUTHOR = {Iskovskih, V. A. and Manin, Ju. I.},
     TITLE = {Three-dimensional quartics and counterexamples to the {L}\"{u}roth
              problem},
   JOURNAL = {Mat. Sb. (N.S.)},
  FJOURNAL = {Matematicheski\u{\i} Sbornik. Novaya Seriya},
    VOLUME = {86(128)},
      YEAR = {1971},
     PAGES = {140--166},
      ISSN = {0368-8666},
   MRCLASS = {14J15},
  MRNUMBER = {291172},
MRREVIEWER = {Peter Br\"{u}ckmann},
}

@Article{Cle77,
 Author = {Clemens, C. H.},
 Title = {Degeneration of {K{\"a}hler} manifolds},
 FJournal = {Duke Mathematical Journal},
 Journal = {Duke Math. J.},
 ISSN = {0012-7094},
 Volume = {44},
 Pages = {215--290},
 Year = {1977},
 Language = {English},
 DOI = {10.1215/S0012-7094-77-04410-6},
 zbMATH = {3549118},
 Zbl = {0353.14005}
}

@incollection {Bea81,
    AUTHOR = {A. Beauville},
     TITLE = {Les singularités du diviseur {$\Theta$} de la jacobienne intermédiaire de l'hypersurface cubique dans {$\mathbb P^4$}},
 BOOKTITLE = {Algebraic threefolds (Varenna, 1981)},
    SERIES = {Lecture Notes in Math.},
    VOLUME = {947},
     PAGES = {190--208},
 PUBLISHER = {Springer, Berlin-New York},
      YEAR = {1982},
}

@article {Bea77,
    AUTHOR = {Beauville, Arnaud},
     TITLE = {Vari\'{e}t\'{e}s de {P}rym et jacobiennes interm\'{e}diaires},
   JOURNAL = {Ann. Sci. \'{E}cole Norm. Sup. (4)},
  FJOURNAL = {Annales Scientifiques de l'\'{E}cole Normale Sup\'{e}rieure. Quatri\`eme
              S\'{e}rie},
    VOLUME = {10},
      YEAR = {1977},
    NUMBER = {3},
     PAGES = {309--391},
      ISSN = {0012-9593},
   MRCLASS = {14K30 (14C30 14J99 14M20)},
  MRNUMBER = {472843},
MRREVIEWER = {T. Oda},
       URL = {http://www.numdam.org/item?id=ASENS_1977_4_10_3_309_0},
}

@article {CG,
    AUTHOR = {Clemens, C. Herbert and Griffiths, Phillip A.},
     TITLE = {The intermediate {J}acobian of the cubic threefold},
   JOURNAL = {Ann. of Math. (2)},
  FJOURNAL = {Annals of Mathematics. Second Series},
    VOLUME = {95},
      YEAR = {1972},
     PAGES = {281--356},
      ISSN = {0003-486X},
   MRCLASS = {14J10 (14G13 14J05 14K20 14M20 14N99)},
  MRNUMBER = {302652},
MRREVIEWER = {H. Popp},
       DOI = {10.2307/1970801},
       URL = {https://doi.org/10.2307/1970801},
}

@article {HT19,
    AUTHOR = {Hassett, Brendan and Tschinkel, Yuri},
     TITLE = {On stable rationality of {F}ano threefolds and del {P}ezzo
              fibrations},
   JOURNAL = {J. Reine Angew. Math.},
  FJOURNAL = {Journal f\"{u}r die Reine und Angewandte Mathematik. [Crelle's
              Journal]},
    VOLUME = {751},
      YEAR = {2019},
     PAGES = {275--287},
      ISSN = {0075-4102},
   MRCLASS = {14E08 (14J30 14J45 14M20)},
  MRNUMBER = {3956696},
MRREVIEWER = {Massimiliano Mella},
       DOI = {10.1515/crelle-2016-0058},
       URL = {https://doi.org/10.1515/crelle-2016-0058},
}

@article {Sho83,
    AUTHOR = {Shokurov, V. V.},
     TITLE = {Prym varieties: theory and applications},
   JOURNAL = {Izv. Akad. Nauk SSSR Ser. Mat.},
  FJOURNAL = {Izvestiya Akademii Nauk SSSR. Seriya Matematicheskaya},
    VOLUME = {47},
      YEAR = {1983},
    NUMBER = {4},
     PAGES = {785--855},
      ISSN = {0373-2436},
   MRCLASS = {14H40 (14K10 14K30)},
  MRNUMBER = {712095},
MRREVIEWER = {Werner Kleinert},
}

@incollection {Tyu80,
    AUTHOR = {Tjurin, A. N.},
     TITLE = {The intermediate {J}acobian of three-dimensional varieties},
 BOOKTITLE = {Current problems in mathematics, {V}ol. 12 ({R}ussian)},
    SERIES = {Itogi Nauki i Tekhniki},
     PAGES = {5--57, 239 (loose errata)},
 PUBLISHER = {Akad. Nauk SSSR, Vsesoyuz. Inst. Nauchn. i Tekhn. Inform.,
              Moscow},
      YEAR = {1979},
   MRCLASS = {14J30},
  MRNUMBER = {537684},
MRREVIEWER = {Miles Reid},
}

@article {Voisin_Quartic7Nodes,
    AUTHOR = {Voisin, Claire},
     TITLE = {Unirational threefolds with no universal codimension {$2$}
              cycle},
   JOURNAL = {Invent. Math.},
  FJOURNAL = {Inventiones Mathematicae},
    VOLUME = {201},
      YEAR = {2015},
    NUMBER = {1},
     PAGES = {207--237},
      ISSN = {0020-9910},
   MRCLASS = {14M20 (14E08 14F45 14J30 14K30)},
  MRNUMBER = {3359052},
MRREVIEWER = {M. Kh. Gizatullin},
       DOI = {10.1007/s00222-014-0551-y},
       URL = {https://doi.org/10.1007/s00222-014-0551-y},
}

@article {Voi88,
    AUTHOR = {Voisin, Claire},
     TITLE = {Sur la jacobienne interm\'{e}diaire du double solide d'indice
              deux},
   JOURNAL = {Duke Math. J.},
  FJOURNAL = {Duke Mathematical Journal},
    VOLUME = {57},
      YEAR = {1988},
    NUMBER = {2},
     PAGES = {629--646},
      ISSN = {0012-7094,1547-7398},
   MRCLASS = {14K30 (14J30)},
  MRNUMBER = {962523},
MRREVIEWER = {A.\ S.\ Tikhomirov},
       DOI = {10.1215/S0012-7094-88-05728-6},
       URL = {https://doi-org.ezproxy.lib.purdue.edu/10.1215/S0012-7094-88-05728-6},
}

@article{Engel2025Matroids,
  author  = {Engel, Philip and Fortman, Olivier de Gaay and Schreieder, Stefan},
  title   = {Matroids and the integral Hodge conjecture for abelian varieties},
  journal = {arXiv preprint arXiv:2507.15704},
  year    = {2025},
  doi     = {10.48550/arXiv.2507.15704}
}

@misc{fanography,
  author = {Belmans, Pieter},
  title = {Fanography},
  howpublished = {\url{https://fanography.info}},
  doi = {10.5281/zenodo.21296718},
  year = {2026},
}

@article {Gwena,
    AUTHOR = {Gwena, Tawanda},
     TITLE = {Degenerations of cubic threefolds and matroids},
   JOURNAL = {Proc. Amer. Math. Soc.},
  FJOURNAL = {Proceedings of the American Mathematical Society},
    VOLUME = {133},
      YEAR = {2005},
    NUMBER = {5},
     PAGES = {1317--1323},
      ISSN = {0002-9939,1088-6826},
   MRCLASS = {14H40 (05B35)},
  MRNUMBER = {2111937},
MRREVIEWER = {Sabin\ Cautis},
       DOI = {10.1090/S0002-9939-04-07689-0},
       URL = {https://doi-org.ezproxy.lib.purdue.edu/10.1090/S0002-9939-04-07689-0},
}

@article {Bardelli,
    AUTHOR = {Bardelli, Fabio},
     TITLE = {Polarized mixed {H}odge structures: on irrationality of
              threefolds via degeneration},
   JOURNAL = {Ann. Mat. Pura Appl. (4)},
  FJOURNAL = {Annali di Matematica Pura ed Applicata. Serie Quarta},
    VOLUME = {137},
      YEAR = {1984},
     PAGES = {287--369},
      ISSN = {0003-4622},
   MRCLASS = {14J30 (14C30 14J15 14K30)},
  MRNUMBER = {772264},
MRREVIEWER = {Steven\ M.\ Zucker},
       DOI = {10.1007/BF01789401},
       URL = {https://doi-org.ezproxy.lib.purdue.edu/10.1007/BF01789401},
}

@article {Hoyt,
    AUTHOR = {Hoyt, William L.},
     TITLE = {On products and algebraic families of jacobian varieties},
   JOURNAL = {Ann. of Math. (2)},
  FJOURNAL = {Annals of Mathematics. Second Series},
    VOLUME = {77},
      YEAR = {1963},
     PAGES = {415--423},
      ISSN = {0003-486X},
   MRCLASS = {14.51 (14.20)},
  MRNUMBER = {150145},
MRREVIEWER = {Y.\ Nakai},
       DOI = {10.2307/1970125},
       URL = {https://doi-org.ezproxy.lib.purdue.edu/10.2307/1970125},
}

@article{Tak89,
  author  = {Takeuchi, Kiyohiko},
  title   = {Some birational maps of {Fano} 3-folds},
  journal = {Compositio Mathematica},
  volume  = {71},
  number  = {3},
  year    = {1989},
  pages   = {265--283},
  url     = {https://www.numdam.org/item/CM_1989__71_3_265_0/}
}

@article{Tyurin,
 author = {Tyurin, A. N.},
 title = {The middle {Jacobian} of three-dimensional varieties},
 fjournal = {Journal of Soviet Mathematics},
 journal = {J. Sov. Math.},
 issn = {0090-4104},
 volume = {13},
 pages = {707--745},
 year = {1980},
 language = {English},
 doi = {10.1007/BF01084562},
 zbMATH = {3666952},
 Zbl = {0428.14015}
}

@Article{AKMW02,
  Author = {Abramovich, Dan and Karu, Kalle and Matsuki, Kenji and W{\l}odarczyk, Jaros{\l}aw},
  Title = {Torification and factorization of birational maps},
  Journal = {J. Amer. Math. Soc.},
  Volume = {15},
  Number = {3},
  Pages = {531--572},
  Year = {2002},
  DOI = {10.1090/S0894-0347-02-00309-2}
}

@Article{Bar84,
  Author = {Bardelli, Fabio},
  Title = {Polarized mixed {Hodge} structures: on irrationality of threefolds via degeneration},
  Journal = {Ann. Mat. Pura Appl. (4)},
  Volume = {137},
  Pages = {287--369},
  Year = {1984}
}

@InCollection{Cle82,
  Author = {Clemens, C. Herbert},
  Title = {Degeneration techniques in the study of threefolds},
  Booktitle = {Algebraic Threefolds (Varenna, 1981)},
  Series = {Lecture Notes in Mathematics},
  Volume = {947},
  Publisher = {Springer},
  Address = {Berlin--New York},
  Pages = {93--154},
  Year = {1982}
}

@Article{Far25,
  Author = {Farb, Benson},
  Title = {Irrationality of the general smooth quartic $3$-fold using intermediate {Jacobians}},
  Journal = {Adv. Math.},
  Volume = {465},
  Pages = {Paper No. 110160, 6},
  Year = {2025},
  DOI = {10.1016/j.aim.2025.110160}
}

@Article{IM00,
  Author = {Iliev, Atanas and Markushevich, Dimitri},
  Title = {The {Abel--Jacobi} map for a cubic threefold and periods of {Fano} threefolds of degree $14$},
  Journal = {Doc. Math.},
  Volume = {5},
  Pages = {23--47},
  Year = {2000}
}

@InCollection{IP99,
  Author = {Iskovskikh, V. A. and Prokhorov, Yu. G.},
  Title = {Fano varieties},
  Booktitle = {Algebraic Geometry V},
  Series = {Encyclopaedia of Mathematical Sciences},
  Volume = {47},
  Publisher = {Springer},
  Address = {Berlin},
  Pages = {1--247},
  Year = {1999}
}

@Article{MM81,
  Author = {Mori, Shigefumi and Mukai, Shigeru},
  Title = {Classification of {Fano} $3$-folds with $B_2\geq 2$},
  Journal = {Manuscripta Math.},
  Volume = {36},
  Pages = {147--162},
  Year = {1981},
  Note = {Erratum: Manuscripta Math. 110 (2003), 407},
  DOI = {10.1007/BF01358754}
}

@Article{Pro18,
  Author = {Prokhorov, Yuri G.},
  Title = {The rationality problem for conic bundles},
  Journal = {Russian Math. Surveys},
  Volume = {73},
  Number = {3},
  Pages = {375--456},
  Year = {2018},
  DOI = {10.1070/RM9811}
}

@misc{Pro26,
  author = {Prokhorov, Yuri G.},
  title = {Double {Veronese} cones with singularities},
  year = {2026},
  eprint = {2605.11932},
  archiveprefix = {arXiv},
  primaryclass = {math.AG},
  note = {\url{https://arxiv.org/abs/2605.11932}}
}

@incollection {ChePar,
    AUTHOR = {Cheltsov, Ivan and Park, Jihun},
     TITLE = {Sextic double solids},
 BOOKTITLE = {Cohomological and geometric approaches to rationality
              problems},
    SERIES = {Progr. Math.},
    VOLUME = {282},
     PAGES = {75--132},
 PUBLISHER = {Birkh\"{a}user Boston, Boston, MA},
      YEAR = {2010},
      ISBN = {978-0-8176-4933-3},
   MRCLASS = {14E08 (14E05 14G05 14G17 14J45)},
  MRNUMBER = {2605166},
MRREVIEWER = {Andreas\ H\"{o}ring},
       DOI = {10.1007/978-0-8176-4934-0\{_}4}

@article {Collino79_irrationality,
    AUTHOR = {Collino, Alberto},
     TITLE = {A cheap proof of the irrationality of most cubic threefolds},
   JOURNAL = {Boll. Un. Mat. Ital. B (5)},
  FJOURNAL = {Unione Matematica Italiana. Bollettino. B. Serie V},
    VOLUME = {16},
      YEAR = {1979},
    NUMBER = {2},
     PAGES = {451--465},
   MRCLASS = {14J30},
  MRNUMBER = {546468},
MRREVIEWER = {Werner\ Kleinert},
}

\end{document}